\documentclass[letterpaper, 10 pt, conference]{ieeeconf} 
\IEEEoverridecommandlockouts                              
\usepackage{graphicx}
\usepackage{epstopdf}
\usepackage{float}
\usepackage{subfig}
\usepackage{color}
\usepackage{todonotes}

\usepackage{amsmath}
\usepackage{amssymb}
\usepackage{mathtools}
\usepackage{algorithm}
\usepackage[noend]{algpseudocode}
\mathtoolsset{showonlyrefs}
\usepackage{textcomp}

\newcommand{\bx}{{\bf x}}

\newcommand{\minimize}{\mathrm{minimize}}
\newcommand{\subject}{\mathrm{subject~to}}

\newtheorem{proposition}{Proposition}
\newtheorem{remark}{Remark}
\newtheorem{theorem}{Theorem}

\usepackage[belowskip=-10pt,aboveskip=5pt]{caption}

\usepackage{tikz}
\usetikzlibrary{shapes,arrows}

\title{\LARGE \bf
Optimal Mass Transport and Kernel Density Estimation \\ for State-Dependent Networked Dynamic Systems
}

\author{Mathias Hudoba de Badyn, Utku Eren, Beh\c{c}et A\c{c}\i kme\c{s}e and Mehran Mesbahi
\thanks{The research of the authors was supported by  NSERC under grant number CGSD2-502554-2017, the U.S. ARL and the U.S. ARO under contract number W911NF-13-1-0340, and the U.S. AFOSR under grant number FA9550-16-1-0022. 
The authors are with the William E. Boeing Department of Aeronautics \& Astronautics, University of Washington, Seattle WA, USA. Emails:
        {\tt\small \{hudomath,ue,behcet,mesbahi\}@uw.edu}~~~~~\textcopyright IEEE 2018}%
}

\begin{document}

\maketitle
\thispagestyle{empty}
\pagestyle{empty}

\begin{abstract}
State-dependent networked dynamic systems are ones where the interconnections between agents change as a function of the states of the agents.
Such systems are highly nonlinear, and a cohesive strategy for their control is lacking in the literature.
In this paper, we present two techniques pertaining to the density control of such systems.
Agent states are initially distributed according to some density, and a feedback law is designed to move the agents to a target density profile.
We use optimal mass transport to design a feedforward control law propelling the agents towards this target density.
Kernel density estimation, with constraints imposed by the state-dependent dynamics, is then used to allow each agent to estimate the local density of the agents.

\end{abstract}

\section{INTRODUCTION}

Networked dynamic systems arise in many synthetic and natural systems in science and engineering~\cite{Mesbahi2010}; in particular multi-agent systems offer an interesting control paradigm.
Each agent augments the system with an additional (local) computational resource, motivating the concept of \emph{distributed controllers \& estimators}~\cite{Acikmese2014,HudobadebadynFillt2017,Li2017}, leading to a notion of \emph{local} control versus \emph{global} control.
The latter seeks to design control laws to guide groups of agents to a desired objective, and the former seeks to design on-board controllers for each agent to facilitate their role in the global control law.

A class of  the more widely used local control laws is called \emph{consensus}, where each agent averages data from their neighbours to compute a parameter related to their objective -- for example, heading, position or a formation center~\cite{Mesbahi2010,Hudobadebadyn2016,Olfati-Saber2004}.
The attractive feature of consensus is how the interconnections between agents -- the so-called \emph{network} or \emph{graph topology} -- affect the performance and agreement characteristics of the algorithm~\cite{Kim2006,Krishnan2017a}.
Graph-theoretic characteristics, such as symmetry, structural balance and graph spectra provide additional insights into the control-theoretic behaviour of consensus ~\cite{Rahmani2009a,Alemzadeh2017,Hudobadebadyn2017,Altafini2013,HudobaDeBadyn2015a}.

Global controllers for multi-agent systems take many forms~\cite{harrison2018optimal}, eg.,
potential field approaches~\cite{Chaimowicz2005}, 
smoothed particle hydrodynamics~\cite{Pimenta2013} and {density control}~\cite{Eren2017}.
In~\cite{Krishnan2017,Krishnan2016}, density control with only relative measurements of position between agents is considered, and 
mean-field control is used to tackle multi-agent interactions by considering a mass flow in the large-$N$ limit~\cite{Albi2016,caines2013mean,Fornasier2013}.

The main interest area of the current paper is networked dynamic systems in which the underlying network is time-varying.
Examples of such systems are switching and proximity-based consensus~\cite{Mesbahi2010,Olfati-Saber2004,Nabi-Abdolyousefi2011} and the Vicsek flocking model~\cite{Tahbaz-Salehi2007}.
State-dependency refers to networks in which the underlying graph varies based on the state of the nodes.
State-dependent networks have been considered in~\cite{Kim2006,Albi2016,Fornasier2013,Mesbahi2002,Mesbahi2005,Awad2015,Awad2018}, but few underlying principles for designing controllers to account for this difficult nonlinearity has been proposed.
A wide range of real-life networked systems are state-dependent.

To tackle this problem, we propose a twofold extension of the work in~\cite{Eren2017}.
First, we consider state-dependent networked dynamic systems, instead of single-integrator dynamics.
Second, we propose a control method for these systems by using a feedback-based density control law that utilizes \emph{optimal mass transport} (OMT).
OMT was considered for linear systems in~\cite{Chen2017}, non-linear systems in~\cite{elamvazhuthi2016optimal} and for density tracking of non-interacting agents in~\cite{Chen2018}.
Our contribution considers OMT for multi-agent systems, in particular  ones with state-dependent dynamics.

In the OMT problem, the initial and final densities $\rho_0$ \& $\rho_1$ of the agents are specified.
The solution to the problem yields a time-dependent density profile with boundary conditions imposed by $\rho_0,\rho_1$, and a velocity field that together satisfy the continuity equation.

We aim to use this velocity field as a feedforward control input to the state-dependent multi-agent system, coupled with a density control feedback law.
Using a modified form of kernel density estimation (KDE) that takes into account the state-dependent dynamics, we will show that the combination of the two control techniques will allow us to propose a physically realizable control strategy for state-dependent networked dynamic systems.

This paper is organized as follows.
The mathematical preliminaries, including notation, optimal mass transport, and density control with the KDE procedure are reviewed in \S\ref{sec:math-prel}.
The problem statement and paper contributions are outlined in \S\ref{sec:probl-stat-contr}.
We present the feedforward controller based on OMT in \S\ref{sec:feedb-contr-state}, and the state-dependent KDE in \S\ref{sec:dens-estim-kern}.
Examples are provided in \S\ref{sec:examples}, and the paper is concluded in \S\ref{sec:conclusion}.

\section{MATHEMATICAL PRELIMINARIES AND BACKGROUND}
\label{sec:math-prel}
\subsection{Mathematical Notation}

We follow the standard graph theory notation listed in~\cite{Mesbahi2010}.
A measure space $(\Sigma,\mathcal{A},\mu)$ is a triple containing a \emph{sample space} $\Sigma\subset \mathbb{R}^n$, a $\sigma$-algebra of subsets of $\Sigma$, and a \emph{measure} $\mu$ that assigns the `size' $\mu(A) \in \mathbb{R}_+$ to a set $A \in \mathcal{A}$.
The \emph{Borel $\sigma$-algebra} $\mathcal{B}$ is generated from the countable unions, intersections and complements of open subsets of $\mathbb{R}^n$.
The \emph{Lebesgue measure} $\lambda$ assigns to a closed interval $[a,b]\subset\mathbb{R}$ the `size' $b-a$; this can be extended to $\mathbb{R}^n$ by considering products of measures.
A measure $\mu$ is called \emph{absolutely continuous} with respect to a measure $\nu$ if $\nu(A) = 0$ $ \implies$ $ \mu(A)=0$ for $A\in\mathcal{A}$; $\mu$ is called \emph{Lebesgue absolutely continuous} if $\nu = \lambda$.
This is denoted $\mu \ll \nu$.
If $\mu \ll \nu$, there exists a function $f:\mathbb{R}\to\mathbb{R}_+$ called a \emph{Radon-Nikodym derivative}, or \emph{density}, of $\mu$ with respect to $\nu$.
It is denoted $f := \frac{d\mu}{d\nu}$, and satisfies the property that for all $A\in\mathcal{A}$, $\mu(A) = \int_A fd\nu$.
Let $\pi(x,y)$ be a joint measure on $\mathcal{X}\times\mathcal{Y}$.
We denote the set of such joint measures $\pi$ by $P(\mathcal{X},\mathcal{Y})$.
The \emph{marginal} $\pi_{x}$ of $\pi$ on $\mathcal{X}$ is defined as the push-forward under the projection map $X$ on $\mathcal{X}$:
$
  \pi_x = X_\# \pi,
$                                                                                     
where  $X(x,y) = x$.
Similarly, the marginal $\pi_y$ of $\pi$ on $\mathcal{Y}$ is given by
$
  \pi_y = Y_\# \pi,
$
where  $Y(x,y) = y$.
We denote the convolution of two functions $f,g$ as $f\star g$, or the convolution of a function $f$ and measure $\mu$ as $f \star \mu$.

\subsection{Optimal Mass Transport}
Informally speaking, the \emph{optimal mass transport} problem is to find a mapping between two densities that minimizes some cost.
Formally speaking, we consider two  measures\footnote{If the measures are Lebesgue absolutely continuous, one can equivalently consider densities $\rho_0,\rho_1$.} $\mu_0,\mu_1$ on $\mathbb{R}^n$ with equal mass: $\int_{\mathbb{R}^n} d\mu_0 = \int_{\mathbb{R}^n}d\mu_1.$
The \emph{optimal mass transport} problem~\cite{Villani2003,Villani2009} is to find a measurable map $ T:\mathbb{R}^n \to \mathbb{R}^n$ taking $\mu_0$ to $\mu_1$ via the following optimization problem:
\begin{align}
\arraycolsep=1.5pt\def\arraystretch{1.2}
  \begin{array}{ll}
    \minimize &   \int_{\mathbb{R}^n} c\left(x,T(x) \right) d\mu_0(x)\\
    \subject &   \int_{x\in A} d\mu_1(x) = \int_{T(x)\in A} d\mu_0(x),
              \forall A \subset \mathbb{R}^n,
  \end{array}\label{eq:16}
\end{align}
where $c$ is a cost function depending on the initial and transported masses.
The constraint in Problem~\eqref{eq:16} means that $\mu_1$ is the \emph{push-forward} measure of $\mu_0$ under the map $T$, in that for each Borel set $B\in \mathcal{B}:= \sigma(\mathbb{R}^n)$, we have that $\mu_1(B) = \mu_0 \left( T^{-1}(B) \right)$.
This is denoted as $  T_{\#} \mu_0 = \mu_1$.

A generalization of Problem~\eqref{eq:16} by Kantorovich is able to pick out the optimal map $T$, if it exists, for a certain class of costs $c$ under the assumption of absolute continuity of the measures~\cite{Kantorovich1948}.
Here, we consider a joint distribution $\pi(x,y)$ on $\mathbb{R}^n\times\mathbb{R}^n$  and solve for the optimal admissible measure $\pi$ given some cost $c(x,y)$.

The set of admissible measures $\pi$ are those whose marginals are $\mu_0$ and $\mu_1$: 
$
  X_\# \pi = \mu_0,~  Y_\# \pi = \mu_1.
$
This is equivalent to requiring 
\begin{align}
  \pi(A \times \mathbb{R}^n) = \mu_0 (A),~\pi(\mathbb{R}^n \times B) = \mu_1(B)\label{eq:17}
\end{align}
for all measurable $A\subset\mathbb{R}^n$ and $B\subset\mathbb{R}^n$.
The \emph{Kantorovich relaxed optimal mass transport problem}~\cite{Kantorovich1948} is given by 
\begin{align}
  \begin{array}{ll}
    \minimize & \int_{\mathbb{R}^n\times \mathbb{R}^n} c(x,y) d\pi(x,y)\\
    \subject & \pi \in \left\{ \pi \in P(\mathbb{R}^n,\mathbb{R}^n)~ \mathrm{s.t.}~\eqref{eq:17}~\text{holds} \right\}.
  \end{array}\label{eq:18}
\end{align}

\begin{proposition}[\cite{Villani2003,Villani2009}]
  For quadratic costs $c(x,y) = \|x-y\|^2$, the support of the optimal joint measure $\pi^*(x,y)$ of Problem~\eqref{eq:18} is exactly the graph of the optimal map $T^*(x)$ minimizing Problem~\eqref{eq:16}.
\end{proposition}

For quadratic costs, Benamou and Brenier formulated an equivalent problem in terms of a constrained fluid mechanics model.
\begin{proposition}[\cite{Benamou2000}]
\label{prop:bb}
  Given Lebesgue absolutely continuous $\mu_0,\mu_1$ with Radon-Nikodym densities $\rho_0,\rho_1$ respectively, Problem \eqref{eq:18} with quadratic costs is equivalent to
  \begin{align}
\arraycolsep=1.6pt\def\arraystretch{1.2}
    \begin{array}{ll}
      \inf_{\rho,v} &\int_{\mathbb{R}^n} \int_0^1 \frac{1}{2} \|v(t,x)\|^2\rho(t,x)dtdx\\
\subject &  \frac{\partial \rho}{\partial t} + \nabla \cdot (v \rho) = 0\\
                    &\rho(0,x) = \rho_0(x),~\rho(1,y) = \rho_1(y).
    \end{array}\label{eq:23}
  \end{align}
Furthermore, the solution to Problem~\eqref{eq:23} is of the form 
$
  v(t,x) = \nabla \varphi(t,x),
$
where $\varphi(t,x)$ is the Lagrange multiplier of the constraints 
and the solution to the Hamilton-Jacobi equation $\partial_t \phi + \frac{1}{2}|\nabla \phi|^2 = 0$.
\end{proposition}
\begin{remark}
  The optimal map $T^*$ of Problem~\eqref{eq:18} in the case of quadratic costs can be reconstructed from the variable $v(t,x)$ from the solution of Problem~\eqref{eq:23}.
This formally establishes the equivalence stated in Proposition~\ref{prop:bb}~\cite{Benamou2000}.
\end{remark}

\subsection{Density Control and Kernel Density Estimation}
In~\cite{Eren2017}, a feedback control law to drive a group of single-integrator agents to a desired density profile $\rho_1(x,t)$ was analyzed.
The following feedback law is proposed to compute the velocity field as a function of the error in density $\Phi(x,t):=\rho(x,t) - \rho_1(x)$ as
\begin{equation}
	v(x,t) = - \hspace{0.5mm} \alpha \hspace{0.5mm} \frac{\nabla\Phi(x,t)}{\rho(x,t)} \hspace{1mm}.
	\label{eqn: controlLaw}
\end{equation} 

Density control of multi-agent systems is impacted by the ability of individual agents to discern the local density profile from measurements of their neighbours.
Since the number of agents is finite, the local density profile must be approximated from finitely many samples $r_i(t)$.

This can be accomplished using \emph{kernel density estimation}~\cite{wand1994kernel,Parzen1962}.
The kernel density estimate $\hat \rho(t,x)$ at any point $x\in\mathbb{R}^n$ and time $t\in\mathbb{R}_+$ is given by 
\begin{align}
  \hat \rho(t,x) = \int_{\mathbb{R}^n} \left[ \prod_{k=1}^d \dfrac{1}{h_k} K \left( \dfrac{x^{[k]} - \xi^{[k]}}{h_k} \right) \right] dP_N(t,\xi). \label{eq:29}
\end{align}
Here, $K:\mathbb{R}\to \mathbb{R}$ is called the \emph{smoothing kernel}, $h_k$ is called the \emph{smoothing parameter}, and $dP_N(t,\xi)$ is a sum of Dirac measures at sample points: 
\begin{align}
  dP_N(t,\xi) = \dfrac{1}{N} \sum_{r(t) \in S(t)} \delta\left(\xi - r(t) \right) d\xi.
\end{align}
Since $\delta(\cdot)$ is the Dirac delta functional, Equation~\eqref{eq:29} can be written as 
\begin{align}
  \hat \rho(t,x) = \dfrac{1}{N} \sum_{i=1}^N \left[ \prod_{k=1}^d \dfrac{1}{h_k} K \left(\dfrac{x^{[k]} - r_i^{[k]}}{h_k} \right) \right].
\end{align}
The control law~\eqref{eqn: controlLaw} then uses the estimate $\hat\rho(x,t)$ in place of knowledge of the true density $\rho(x,t)$, where the sample points $r(t)$ are taken to be the agent states:
\begin{align}
  v(x,t) = - \alpha \dfrac{\nabla( \hat \rho(x,t) - \rho_1(x) )}{\hat \rho (x,t) }.
\end{align}

In~\cite{Eren2017}, Gaussian kernels were used -- this induces an all-to-all communication; every agent is able to sample the position every other agent.
The control law~\eqref{eqn: controlLaw} has a  convergence guarantee listed in Theorem~6 of \cite{Eren2017}.

\section{PROBLEM STATEMENT AND CONTRIBUTIONS}
\label{sec:probl-stat-contr}

In this paper, we consider state-dependent networked dynamic systems with $N$ agents on a bounded region $\mathcal{R}\subset \mathbb{R}^n$, where the $i$th agent's state evolves according to the dynamics
\begin{align}
  \dot x_i = f(\mathcal{G}(\textbf{x}),\textbf{x}) + B u_i(\textbf{x},t),~\textbf{x}:=(x_1,\dots,x_n).
\end{align}
A prototypical example of such a system is \emph{state-dependent consensus} 
\begin{align}
  \dot x_i =\sum_{j\neq i} A(x_i,x_j)\cdot(x_i - x_j) + u_i,\label{eq:2}
\end{align}
where the edge weight $w_{ij}:=A(x_i,x_j) $ changes depending on the state of the agent $i$ and its neighbour $j$.
In general, one can consider an \emph{interaction kernel} $H(x)$ that generates a consensus-like dynamics by convolution with the Dirac measure supported at agent states~\cite{Albi2016,Fornasier2013}: 
\begin{align}
  \mu_N(x) = \frac{1}{N}\sum_{j=1}^N \mathbf{1}_{\{x_j\}}(x) = \frac{1}{N} \sum_{j=1}^N \delta(x-x_j).\label{eq:1}
\end{align}
Using Equation~\eqref{eq:1}, we can write a general multi-agent system as 
\begin{align}
  \dot x_i = \left(H \star \mu_N\right)(x_i) + u_i.
\end{align}
A simple example motivated by robotics is proximity-based edge switching, where $A(x_i,x_j) = 1$ if $\|x_i - x_j \| \leq r$ and 0 otherwise, where $r$ is some communication radius. 
This corresponds to an interaction kernel 
$
 H(x) = x \mathbf{1}_{\|x\|\leq r} (x).
$

As the number of agents $N$ grows sufficiently large, one can consider the time-dependent \emph{density} of agents $\rho(x,t)$ over a region of the state space.
In the context of mean-field control, the formal large-$N$ limit of the dynamics~\eqref{eq:2} produces the \emph{mean field dynamics} \cite{Albi2016,Fornasier2013},
\begin{align}
\arraycolsep=0.5pt\def\arraystretch{1.5}
  \begin{array}{ll}
    &  \dfrac{\partial \rho }{\partial t} + \nabla_x\cdot \left[ \left( \mathcal{P}(\rho(x,t),t) + u\right) \rho \right] = 0\\ 
& \mathcal{P}(\rho(x,t),x) = \int A(x,y)(y-x) \rho(y,t)dy.
  \end{array}\label{eq:3}
\end{align}
We now state the contributions of this paper, namely the feedforward OMT scheme with kernel density estimation shown in Figure~\ref{fig:blkdiag}.

\tikzstyle{block} = [draw, fill=blue!20, rectangle, 
    minimum height=2em, minimum width=4em]
\tikzstyle{sum} = [draw, fill=blue!20, circle, node distance=0.75cm]
\tikzstyle{input} = [coordinate]
\tikzstyle{output} = [coordinate]
\tikzstyle{pinstyle} = [pin edge={to-,thin,black}]
\begin{figure}
  \centering
  \begin{tikzpicture}[auto, node distance=2cm,>=latex']
    \node [input, name=input] {};
    \node [block, right of=input] (omt) {OMT};
    \node [sum, right of=omt,node distance=2cm] (sum) {};
    \node [block, right of=sum] (swarm) {Swarm};
    \draw [->] (omt) -- node[name=u] {$\rho(x,t)$} (sum);
    \node [output, right of=swarm] (output) {};
    \node [block, below of=swarm,node distance=1cm] (measurements) {KDE};

    \draw [->] (sum) -- node {$v(x,t)$} (swarm);
    \draw [->] (swarm) -- node [name=y] {$\bf{x}$}(output);
    \draw [->] (y) |- (measurements);
    \draw [->] (measurements) -| node[pos=0.4] {} node [pos=0.25, above] {$\hat{\rho}(x,t)$} (sum);
      \end{tikzpicture}
\caption{Block diagram of density control scheme}
\label{fig:blkdiag}
\end{figure}
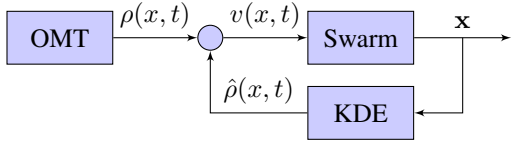

\textit{Contribution 1: Feed-Forward Density Control With Optimal Mass Transport.}

We consider a generalization of the Brenier-Benamou OMT problem \eqref{eq:23} with the continuity equation constraint replaced by the mean-field dynamics~\eqref{eq:3}:
\begin{align}
  \arraycolsep=1.6pt\def\arraystretch{1.2}
  \begin{array}{ll}
    \inf_{\rho,v} &\int  \int_0^1 \frac{1}{2} \|v(t,x)\|^2\rho(t,x)dtdx\\
    \subject    &  {\partial_t \rho } + \nabla_x \cdot \left[ \left( \mathcal{P}(\rho(x, t),x) + v\right) \rho \right] = 0\\ 
                  & \mathcal{P}(\rho(x,t),x) = \int A(x,y)(y-x) \rho(y,t)dy\\
                  &\rho(0,x) = \rho_0(x),~\rho(1,y) = \rho_1(y).
  \end{array}\label{eq:4}
\end{align}
The solution to Problem~\eqref{eq:4} yields two variables with important physical interpretations.
The time-varying density $\rho(t,x)$ represents the mass of agents with dynamics \eqref{eq:2} with inputs $u_i(x,t):= v(x_i,t)$; the velocity field $v(x,t)$ is precisely the input $u_i$ given to agent $i$ at position $x$ at time $t$.

Problem~\eqref{eq:4} assumes that the initial and final masses of agents are distributed in the mean-field limit according to densities $\rho_0$ and $\rho_1$, respectively.
When considering finitely many agents, any initial density will take the form of Equation~\eqref{eq:1} - namely, it will be a Dirac measure supported at the agent states, also called the \emph{empirical density}.

Hence, in general, the boundary conditions on the density in Problem~\eqref{eq:4} can be either \emph{deterministic} (in the case of Dirac measures supported at the agent states), or \emph{probabalistic} (in the sense that the initial/final agent states $x_i(0)$ and $x_i(1)$ are randomly distributed according to the densities $\rho_0$ and $\rho_1$).
In the latter case, as $N\to\infty$, the Dirac measure supported at the agent states at time $t$ converges in a formal sense to the density $\rho(\cdot,t)$~\cite{Albi2016,Fornasier2013}.

In either case, we consider the velocity field $v(x,t)$ as a feed-forward input to the dynamics~\eqref{eq:2}.
Since the number of agents is finite, the empirical density at time $t$ will only approximate the density $\rho(x,t)$ from the solution of Problem~\eqref{eq:4}.

\textit{Contribution 2: Feedback Density Control with Kernel Density Estimation and State-Dependent Constraints.}

In a state-dependent networked dynamic system, eg., Equation~\eqref{eq:2}, the existence of an edge indicates some notion of information transfer between agents.
Hence, a physical estimation scheme and density control law can only allow $i$ to sample those agents $j$ such that $A(x_i,x_j) \neq 0$.
The second contribution of this paper is to extend the KDE procedure in~\cite{Eren2017} by solving a quadratic program for an optimal kernel that takes into account the state-dependent communication constraints.

\section{FEEDBACK CONTROL OF STATE-DEPENDENT NETWORKED DYNAMIC SYSTEMS}
\label{sec:feedb-contr-state}

Consider the state-dependent consensus dynamics \eqref{eq:2}.
Let $v_1(x,t)$ denote the velocity field from the solution to Problem~\eqref{eq:4}, and let $v_2(x,t)$ denote the velocity field from the control law~\eqref{eqn: controlLaw}.
Our proposed control law is then given by the velocity field (with $\alpha > 0$),
\begin{align}
  &u(x,t) = 
           \begin{cases}
             v_1(x,t) - \alpha \dfrac{\nabla\left( \rho(x,t) - \rho_1(x) \right)}{\rho(x,t)} &  0 \leq t \leq  1\vspace{1mm}\\
             v_2(x,t) - \mathcal{P}(\mu(x,t),x) & t \geq 1
           \end{cases}\label{eq:5}\\
&\mathcal{P}(\mu(x,t),x) = \int A(x,y)(y-x) \mu(y,t)dy.
\end{align}
Of course, the switch at $t=1$ is completely arbitrary, and can be altered by changing the time horizon of Problem~\eqref{eq:4}.

The main result of this section is the following theorem, an extension of Theorem 6 in~\cite{Eren2017}.
Informally, it states that as the number of agents in the system tends as $N\to\infty$, the velocities of the agents performing the state-dependent control law~\eqref{eq:5} will vanish asymptotically.
\begin{theorem}
  \label{thr:1}
  \textit{
  Consider a system of $N$ agents $\mathcal{S}(t)$ on a bounded region $\mathcal{R}\subset\mathbb{R}^n$ with individual dynamics given by \eqref{eq:2} and with control law~\eqref{eq:5}.
  Further suppose that the initial swarm density $\rho(x,t)$ and target density $\rho_1(x)$ satisfy the boundary condition $\nabla \Phi(x,t) = 0$ on $\partial \mathcal{R}$.
  As $t\to\infty$, for sufficiently large $N$ the error density $\Phi(x,t)$ converges to zero: 
$  \lim_{t\to\infty} \Phi(x,t) = 0,~\text{for }x \in \mathcal{R}$
and so $\hat\rho(x,t) \to \rho_1(x)$.
Furthermore, the velocities of all agents vanish asymptotically: 
$
  \lim_{t\to\infty} \dot x(t) = 0, ~\text{for }x\in \mathcal{R}.
$
}
\end{theorem}
The proof is discussed in the Appendix.

\section{DENSITY ESTIMATION FOR KERNELS\\ WITH COMPACT SUPPORT}
\label{sec:dens-estim-kern}

Consider a state-dependent consensus dynamics as in Equation~\eqref{eq:2}
where $A(x_i,x_j)$ is a state-dependent edge weight.
In order to implement the density control law, each agent must be able to estimate the density of nearby agents to generate the correct velocity field.
In this section, we discuss optimal kernels designed to achieve this task that are subject to the state-dependent constraints imposed by $A(x_i,x_j)$.

\begin{figure}
  \centering
  \includegraphics[scale = 0.5]{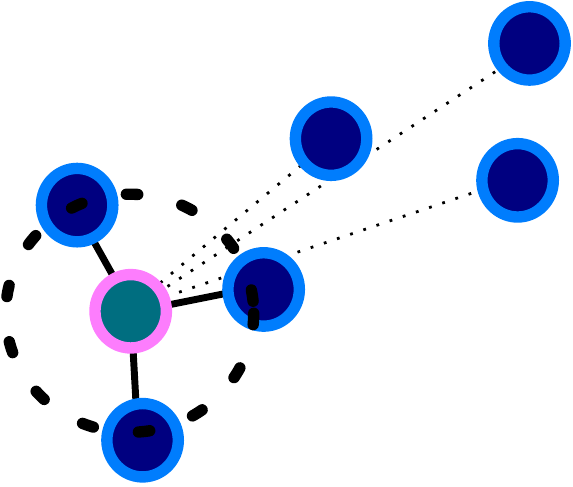}
  \caption{Illustration of (proximity-based) state-dependent constraints on the KDE procedure. Dotted lines indicate samples the center agent cannot measure.}
  \label{fig:sub1}
\end{figure}

The state-dependent constraints in some (informal) sense denote `information transfer' between agents $i$ and $j$.
If $A_{ij}=0$, then agents $i$ and $j$ cannot detect each other, and the KDE procedure should reflect this.
To illustrate this notion, consider Figure~\ref{fig:sub1}.

In 1D kernel density estimation, there are two parameters selected \emph{a priori} that influence the quality of the estimated probability density function, namely the kernel $K$ and the smoothing parameter $h$.
We consider an optimal selection of $K$ subject to the state-dependent constraints; we leave the task of selecting $h$ for future work.
For now, we just need the following assumption on $h$ as a function of the number of samples: 
$
  \lim_{N\to\infty} Nh(N) = \infty.
$
The standard metric for measuring the quality of the estimated probability density function is given by the \emph{mean integrated square error} 
$
 E_{\text{MISE}}:=  \mathbb{E} \left[ \int \left( \hat \rho(x) - \rho(x) \right)dx \right].
$

By extracting out the dependence on the number of samples $N$, and the choice of smoothing parameter $h$, one can obtain the \emph{asymptotic mean integrated square error} (AMISE) \cite{wand1994kernel}:
$
  E_{\text{MISE}} := E_{\text{AMISE}} + o \left( (hn)^{-1} + h^4 \right).
$
One can factor the AMISE into a product of two terms, one depending on $h$ and one depending on $K$: $  E_{\text{AMISE}} = C_1(K)C_2(h)$, where
\begin{align}
  C_1(K) &:= \left[ \left(\int K(x)^2 dx \right)^4 \left( \int x^2 K(x) dx \right) ^2\right]^{1/5}.
\end{align}

Hence, by fixing $a:= \int x^2 K(x)dx$ depending on the length of the boundary of our estimation horizon, the only parameter left to optimize is the \emph{roughness} of $K(x)$: $\int K(x)^2dx$.
We can write an optimization problem as follows:
\begin{align}  
  \begin{array}{ll}
    \minimize & \int K(x)^2 dx \\
    \subject & \int K(x) = 1,~ \int xK(x) = 0\\
              & \int x^2 K(x) = a^2 < \infty,~ K(x) \geq 0.
  \end{array}\label{eq:37}
\end{align}
In one dimension, the solution to Problem~\eqref{eq:37} is given by \cite{wand1994kernel},
\begin{align}
  K^a(x) = \dfrac{3}{4} \dfrac{1}{a\sqrt{5}} \left[ 1- \left( \dfrac{x}{a\sqrt{5}} \right)^2 \right] \mathbf{1}_{\{|x| \leq a\sqrt{5}\}}.\label{eq:39}
\end{align}
We wish to find the solution to a modified version of this problem where we enforce a compact support constraint of the form $\{K(x) = 0,~x\in \mathcal{A},~\mathcal{A}^c~\text{compact}\}$.
To this end, notice that  Problem~\eqref{eq:37} can be numerically solved by discretizing it as follows.

Denote the region of the problem as $X= \{ x ~:~ |x| \leq B\}$.
Discretize $X$ into $N$ points spaced $dx$ apart.
Let the vector $x:=\{x_i\}_{i=1}^N \in [-B,B]^N$ consist of these points, and let $k:=\{k_i\}_{i=1}^N$ be the vector of the kernel $K$ evaluated at these points, i.e. $k_i = K(x_i)$.
Discretizing the integrals  yields a quadratic program of the form:
\begin{align}  
  \begin{array}{ll}
    \minimize & k^Tk\\
    \subject & \sum_{i=1}^N k_i dx = 1,~ \sum_{i=1}^N x_i k_i = 0\\
              & \sum_{i=1}^N x_i^2 k_i dx = a^2,~ k_i \geq 0 ,~1\leq i \leq N.
  \end{array}\label{eq:38}
\end{align}

\begin{figure}
  \centering
  \includegraphics[width=\columnwidth]{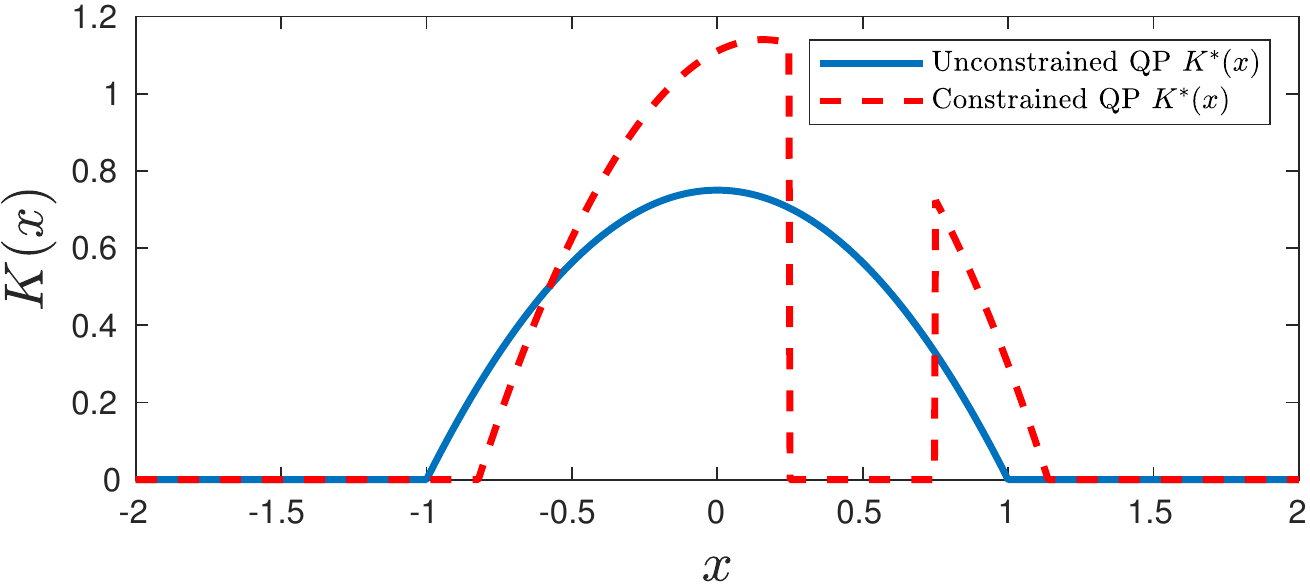}
  \caption{Optimal kernels with unconstrained support, and support constrained to $[-2,2]\setminus{[1/4,3/4]}$, with second moment $a=5^{-1/2}$. The unconstrained kernel solution is exactly given by Equation~\eqref{eq:39}.}
  \label{fig:kernel}
\end{figure}

As discussed before, an agent's state-dependent density estimate will depend on sampling points from agents that have an edge between them.
Hence, in the density estimate for agent $i$, the kernel $K$ will only depend upon the state of agent $i$ and its neighbours $\mathcal{N}_i$.
The density estimate of agent $i$ is written as 
\begin{align}
  \hat{\rho}_i(t,x(t)) = \dfrac{1}{Nh^d} \sum_{j \in\mathcal{N}_i} \left [ \prod_{k=1}^d K_i\left( \dfrac{x_i(t) - x_j(t)}{h}, x_j\right) \right]
\end{align}
where the support of the kernel is restricted to the support of the state-dependent edge weight $A(x_i,x_j)$: 
\begin{align}
 A(x_i,x_j) = 0 \implies K(h^{-1}(x_i-x_j), x_j) = 0. 
\end{align}

To extend Problem~\eqref{eq:38} to multi-dimensional systems, we consider  \emph{multiplicative kernels} for $x_i \in \mathbb{R}^n$, where each dimension is estimated independently:
$
  K(\bx) = \prod_{k=1}^N K_k(x_k).
$
This yields the final optimal kernel problem with compact support constraint.
For brevity, we show the explicit form for the 2D problem, as it is clear (yet notationally cumbersome) how to write the general $N$D problem:
\begin{align}  
  \begin{array}{ll}
    \minimize & \sum_i^{N_x} \sum _j^{N_y} k_{ij}^2\\
    \subject  & \sum_{i=1}^{N_x}\sum_{j=1}^{N_y} k_{ij} dxdy = 1,~k_{ij} \geq 0 ,~\forall i,j\\
              & \sum_{i=1}^{N_x} x_i k_{ij} = 0,~\sum_{i=1}^{N_x} x_i^2 k_{ij} dx = a^2,~\forall j\\
              & \sum_{j=1}^{N_y} y_j k_{ij} = 0,~\sum_{j=1}^{N_y} y_j^2 k_{ij} dx = a^2,~\forall i\\
              &k_{ij} = 0 \text{ if } (x_i,y_j) \in \mathcal{A},~\mathcal{A}^c\text{ compact} .
  \end{array}\label{eq:7}
\end{align}
It is important to note that removing a compact interval from the kernel may bias the density estimate - this is unavoidable.
The compact support constraint defines a \emph{selection-biased distribution (SBD)}; each agent samples the distribution $g(x) =  w(x) \rho(x)/\mu$, with $ w(x) = \mathbf{1}_{\mathcal{A}^c}(x), \mu =\int w(x) \rho(x) dx $.
The standard \emph{unbiased} kernel density estimate of a SBD involves multiplying the kernel by a factor of $\mu/w(x)$~\cite{Gill1988,Wu1996}, which is unbounded for our choice of $w(x)$.
Techniques for unbiasing $\hat\rho(x)$ will be left for future work.


\section{EXAMPLES}
\label{sec:examples}

We numerically simulate $N=200$ agents with interaction kernel $H(x) = x \mathbf{1}_{\|x\|\leq 0.01} (x)$.
The velocity field $\mathcal{P}[\mu](x)$ is evaluated with the Dirac measure \eqref{eq:1}, effectively yielding $N$ single-integrator agents that are able to only sample agents a short distance away from each other.

\begin{figure}
  \centering
  \includegraphics[width=0.6\columnwidth]{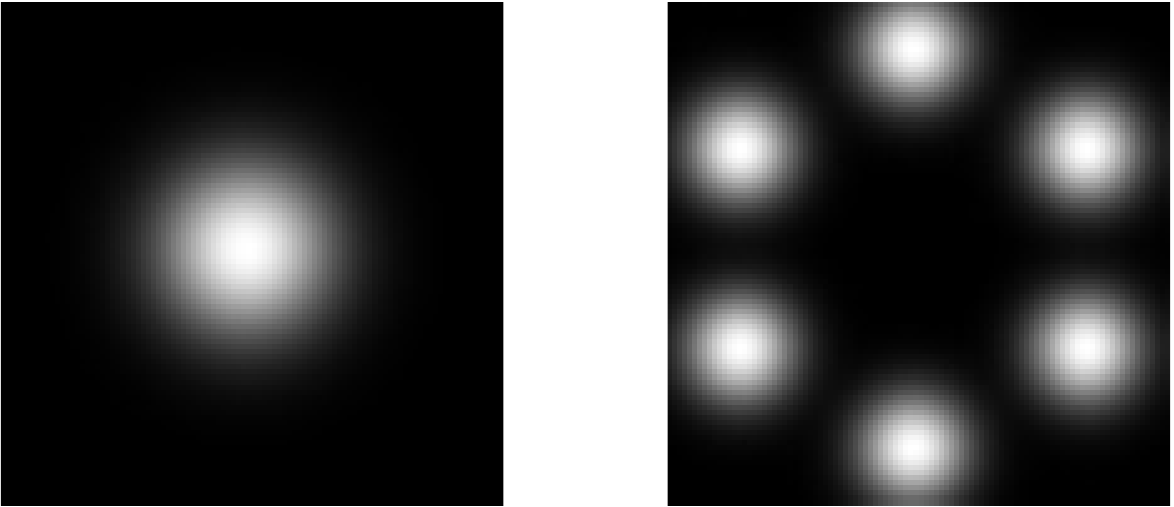}
  \caption{Left: Initial density $\rho_0$. Right: Target density $\rho_1$.}
  \label{fig:dens}
\end{figure}

\begin{figure}
  \centering
  \includegraphics[width=0.65\columnwidth]{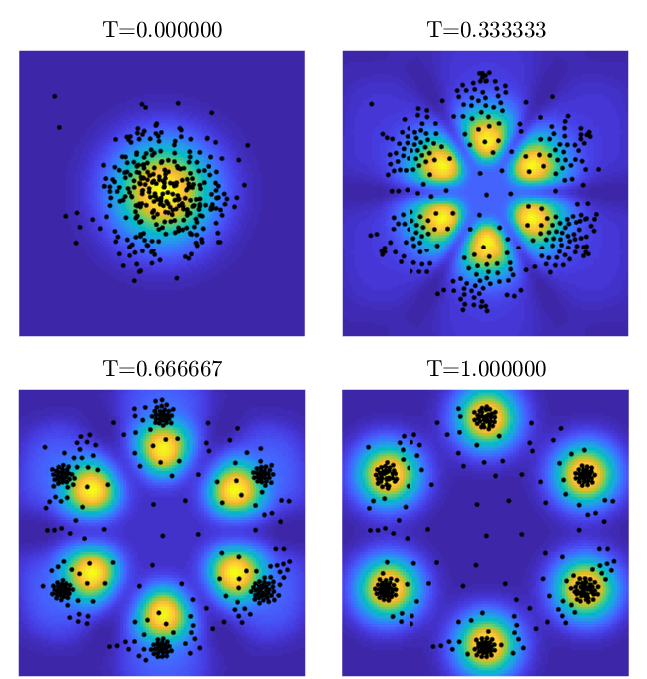}
  \caption{Optimal density profiles over time, and superimposed agent states using the feedback density control law.}
  \label{fig:results}
\end{figure}

\begin{figure}
  \centering
  \includegraphics[width=0.65\columnwidth]{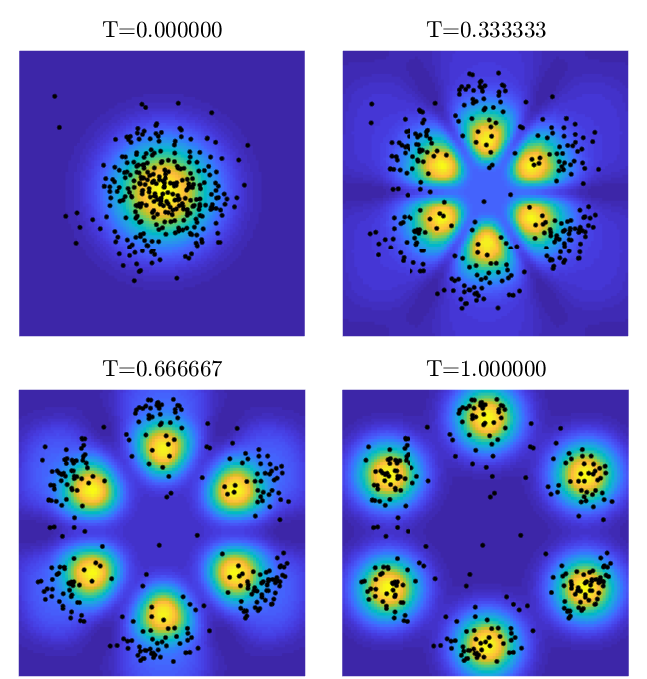}
  \caption{Optimal density profiles over time, and superimposed agent states with only the feedforward control.}
  \label{fig:results2}
\end{figure}

The optimal mass transport problem was solved using open-source code, utilizing a primal-dual algorithm~\cite{Papadakis2014,Peyre2011}.
The density profile and velocity field was calculated over a $100\times100\times100$ grid in $x,y$ and $t$ space.
The initial density $\rho_0(x)$ was a 2D Gaussian at the center of a $[0,1]^2$ grid, and the target density was a ring of 2D Gaussians, as shown in Figure~\ref{fig:dens}.
The superimposed optimal density profile over time, and the states of the agents over time integrating the feedback and feedforward control law are shown in Figure~\ref{fig:results}.
As one can see, the agents are more organized around the final density distributions when using the feedback law as opposed to just integrating the feedforward law, as shown in Figure~\ref{fig:results2}.

\section{CONCLUSION}
\label{sec:conclusion}

In this paper, we examined density control of state-dependent networked dynamic systems.
We utilized the optimal mass transport problem to design a feed-forward velocity field propelling agents with initial conditions sampled from a density profile $\rho_0$ to some target density $\rho_1$.
We then tackled the problem of using a density feedback control law with sparse measurements dictated by the state-dependent edge switching constraints of the agents.
We utilized kernel density estimation to convert measurements of neighboring agents into a local estimate of the swarm density, which was then used to calculate a feedback density control law.
In particular, a quadratic program was designed to find the optimal kernel subject to the state-dependent edge switching.

There are many open problems remaining, here we discuss several.
First, the selection of an optimal interaction distance $r = h$ for proximity-based edge switching. 
This will depend on, for example, $\rho_0,~\rho_1$ and $N$.
If $h$ is large, this will require more on-board computation and sensing capability; if $h$ is small, agents will be isolated.
Second, one can consider the task of determining a state-dependent kernel yielding an unbiased estimate of $\rho$.

\section*{APPENDIX}
We now state the proof of Theorem~\ref{thr:1}.

\begin{proof}
  First, recall the following theorem about consistency of the estimate $\hat \rho(x,t)$. 
  \begin{theorem}[\cite{Wied2012}]
    \label{thr:2}
    \textit{
    Consider a kernel density estimation scheme for the target density $\rho(x)$.
    Suppose the smoothing parameter $h$ is chosen as a function of the number of samples $N$: 
$
      \lim_{N\to\infty} Nh(N) = \infty.
$
    Then, at each point of continuity $x$ of $\rho$, the estimator $\hat \rho_N(x)$ is weakly consistent in that for all $\epsilon>0$, 
$
      \lim_{N\to\infty} \mathbb{P} \left(|\rho_N(x) - \rho(x) | > \epsilon \right) = 0.
$
}
  \end{theorem}
By Theorem~\ref{thr:2}, under the assumption on the smoothing parameter $h$, we have that as $N\to\infty$, $\hat\rho(t,x) \to \rho(t,x)$ with probability 1 for any finite $t$.

Consider the following Lyapunov function:
\begin{align}
  V(t) = \int_{\mathcal{R}} \left( \dfrac{\rho(x,t)}{\alpha} \right)^2 \dot x ^T \dot x~ dx.\label{eq:6}
\end{align}
As $N\to\infty$, the velocity field $\dot x$ approaches the mean-field limit \cite{Albi2016}:
$
  \dot x = \mathcal{P}(\mu,t) + u(x,t),
$
where 
$
  \mathcal{P}(\mu,x) = \int A(x,y)(y-x)\mu(y,t)dy,
$
and $\mu:=\mu(y,t)$ is the measure satisfying 
$
{\partial t} + \nabla \cdot \left( \left( \mathcal{P}(\mu,t) + u(x,t) \right) \mu \right) = 0.
$
Under the control law~\eqref{eq:5}, it follows that for sufficiently large $t$, the Lyapunov function \eqref{eq:6} can be written as 
\begin{align}
  V(t) = \int_{\mathcal{R}} \nabla \Phi(x,t)^T \nabla \Phi(x,t) dx.
\end{align}
The time derivative of $V(t)$ is then given by 
\begin{align}
  \dot V(t)
            &= \alpha \int_{\mathcal{R}} \xi(x,t)^T \Delta \xi(x,t) dx,
\end{align}
where $\xi(x,t) := \nabla \Phi(x,t)$.
Since $\nabla\Phi(x,t) = 0$ on $\partial\mathcal{R}$, we have that $\xi(x,t) = 0 $ on $\partial \mathcal{R}$ which is a Dirichlet boundary condition.
It follows that $\dot V(t) < 0$ since the Dirichlet problem for the Laplace operator has strictly negative eigenvalues~\cite{Li1983}.

Therefore, by LaSalle's Invariance Principle, we can conclude that $\lim_{t\to\infty} \nabla \Phi(x,t)=0$, and so $\lim_{t\to\infty} \Phi(x,t) = \text{constant}$.
However, since $\nabla\Phi(x,t)=0$ on $\partial\mathcal{R}$, the mass $\int_{\mathcal{R}} \Phi(x,t) dx$ is conserved for all $t>0$ (in that $\int_{\mathcal{R}} \Phi(x,t) dx = 0$) and so we have that 
$
  \int_{\mathcal{R}} \hat \rho(x,t) dx = \int_{\mathcal{R}} \rho(x,t) dx.
$
Consequently, it follows that $\lim_{t\to\infty} \Phi(x,t) = 0$, and hence $\lim_{t\to\infty} \rho(x,t) = \lim_{t\to\infty} \hat \rho(x,t)$ which completes the proof.
\end{proof}


\end{document}